\documentclass{amsart}
\usepackage{amssymb}

\newcommand{\bai}{bounded approximate identity\;}

\newcommand{\lcg}{locally compact group\;}

\newcommand{\K}{K\backslash G/H}

\newcommand{\ri}{\rightarrow}

\newcommand{\su}{\subseteq}

\newcommand{\nhd}{neighbourhood\;}
\newcommand{\cpt}{compact\;}

\newcommand{\lcpt}{locally compact\;}

\newcommand{\fn}{function\;}
\newcommand{\fns}{functions\;}

\newcommand{\repn}{representation\;}
\newcommand{\repns}{representations\;}

\newtheorem{defi}{Definition}[section]
\newtheorem{prop}{Proposition}[section]
\newtheorem{theo}{Theorem}[section]
\newtheorem{lemm}{Lemma}[section]
\newtheorem{cor}{Corollary}[section]

\newtheorem{rem}{Remark}[section]

\begin{document}

\title[double coset spaces]{Harmonic analysis on double coset spaces}
\date{}
\author[M. Amini, A.R. Medghalchi]{Massoud Amini, Alireza Medghalchi}
\address{Department of Mathematics\\Shahid Beheshti University\\ Evin, Tehran 19839, Iran\\ m-amini@cc.sbu.ac.ir
\linebreak
\indent Department of Mathematics and Statistics\\ University of Saskatchewan
\\ 106 Wiggins Road\\ Saskatoon, Saskatchewan\\
Canada S7N 5E6\\ mamini@math.usask.ca
\linebreak
\linebreak
\indent Department of Mathematics\\ 
Teacher Training University, Tehran, Iran\\ 
\linebreak
medghalchi@saba.tmu.ac.ir}
\thanks{This research was supported by Grant 510-2090
of Shahid Behshti University}

\begin{abstract}
We prove the existance and uniqueness of quasi-invariant measure on double cost space $\K$ and study the Fourier and Fourier-Stieltjes algebras of these spaces

\end{abstract}
\maketitle

\section{introduction}

Let $G$ be a locally compact group and $H,K$ be closed subgroups
of $G$. The {\it double coset space} of $G$ by $H$ and $K$ is
$$
\K =\{KxH: x\in G\}.
$$
When $K=\{1\}$, this is the {\it coset space} $G/H$ and when $K=H$, it
is the double coset space $G//H$ of $G$ by $H$. Although these are
not groups when $H$ and $K$ are not normal, but a major part of
the Harmonic analysis on $G$ carries over these spaces. Among
these the cost space $G/H$ is quite well studied by several
authors (see for instance [Fl], [RS], [Fo], [Sk] and references
therein). Some part of these works has been carried over the
double coset spaces [Li]. One important point, which is partly
overlooked, is the fact that the double coset space $G//H$ is a
hypergroup (convo) and the coset space $G/H$ is a semi hypergroup
(semi convo) [Je]. 

In this paper we study double coset spaces.
In particular we are interested in studying the function algebras
on these spaces. One major motivation for this study is the fact
that for a general hypergroup, the set of positive definite
functions is not closed under pointwise multiplication [Vr]. In
particular the Fourier and Fourier-Stieltjes spaces fail to be
Banach algebras in general. Since there is no involution on $G/H$,
the concept of positive definiteness could not be defined directly
on this semi hypergroup (there are some alternative definitions,
see for instance [Ar], [Fo]). The double coset space has neither
of these pathologies. It has a involution and positive definite
functions could be defined on it and they have the natural
relation with their counterparts on $G$ [Je]. 
In section 4 we study the Fourier and
Fourier-Stieltjes spaces over $G//H$. Also along the line of
[Fo], we consider the corresponding subalgebras of the Fourier and
Fourier-Stieltjes algebras over $G$ consisting of functions which
are constant on double cosets of $H$.

The paper is organized as follows. In section 2 we study the
(quasi) invariant measures on $\K$ and the analog of the Weil's
formula. The results of this section are mainly contained in [Li],
but as our results are slightly more general and our approach is
different at some points, and also because we later use some of
the technical lemmas in this section which are not included in
[Li], we preferred to present the details of proofs. Section 3 is
a brief account of the theory of induced \repns (Mackey machine)
for $\K$. The first part includes new results along the line of
the classical theory for $G/H$. The second part considers the
induction procedure from the subhypergroups of $G//H$. The results
of this part are taken from [HKK] and [He], where these results
are given in a more general setting. We only quote those results
which are used in the next section. In section 4 we study function
spaces on $\K$ and $G//H$. In particular we are interested in
different versions of Banach spaces defined using positive
definite functions. The proofs in this section partly follow the
same line of reasoning as their coset space counterparts in [Fo].
But there are technical difficulties which force us to present the
details of the proofs and show the spots where the arguments break
down. There is however a major difference. In contrast with the
case of $G/H$, we have an involution structure on $G//H$ and so we
are able to define the function spaces directly on $G//H$, where
as in [Fo], the function (and operator) algebras on $G/H$ are
defined indirectly using the corresponding appropriate subalgebras
of the same algebras on $G$. Therefore from the hypergroup point
of view, our study gives more direct information about the
structure of $G//H$.

Our basic references for the harmonic analysis on $G$ are [Fl] and
[Ey] and we use their,now classical, notations and results without
reference.

\section{quasi-invariant measures}

Let $G$ be \lcg and $H$ ,$K$ be closed subgroups of $G$ and  $dx$, $dh$, and $dk$ denote the Haar measures on them, respectively. If either of these groups are compact we assume that the corresponding Haar measure is normalized so that the measure of the group is one. Let $\Delta_G$, $\Delta_H$, and  $\Delta_K$ denote the corresponding modular functions. The {\it double quotient space of} $G$ {\it by} $H$ {\it and} $K$ is denoted by $\K$ and is defined by
$$
\K=\{KxH:x\in G\}.
$$
Consider the map $q:G\rightarrow\K$, defined by $q(x)=KxH=:\dot{x}\quad (x\in G)$, then the quotient topology of $\K$ is the weakest topology on $\K$ which makes $q$ continuous. This makes $\K$ a \lcpt space on which $G$ acts by translations. For $f\in C_c(G)$ define
$$
Qf(KxH)=\int_K\int_H f(kxh)dkdh,
$$
then $Q:C_c(G)\ri C_c(\K)$ and $supp(Qf)\subseteq q(supp(f)) \quad (f\in C_c(G))$. Note that, by Fubini's theorem, 
the order of integration in above formula is not important. The following lemma is trivial.
\begin{lemm}
For each $f\in C_c(G)$ and $\varphi\in C_c(\K)$
$$
Q((\varphi\circ q).f)=\varphi .Qf.
$$
\qed
\end{lemm}

\begin{lemm}
If $E\su \K$ is \cpt then there is a \cpt set $F\su G$ with $q(F)=E$.
\end{lemm}
{\bf Proof} Let $V$ be a \nhd of identity in $G$ with compact closure and cover $E$ by the sets $q(xV)$, $x\in G$. Choose a subcover $\{q(x_jV)\}^n _1$ and put $F=q^{-1}(E)\cap(\cup_{j=1}^n x_j\bar{V})$.\qed

\begin{lemm}
If $F\su \K$ is \cpt then there is $f\in C_c(G)_{+}$ such that $Qf=1$ on $F$.
\end{lemm}
{\bf Proof} Let $E$ be a \cpt \nhd of $F$ in $\K$ and choose by above lemma a \cpt set $F_1$ in $G$ such that $q(F_1)=E$. Choose $g\in C_c(G)_{+}$ such that $g>0$ on $F_1$ and using Uryshon's Lemma, choose $\varphi\in C_c(\K)$ with $supp(\varphi)\su E$ and $\varphi =1$ on $F$. Put $f=\frac{\varphi\circ q}{Qg\circ q}.g$, then since $Qg>0$ on $supp(\varphi)$, we have $f\in C_c(G)$ and clearly $supp(f)\su supp(g)$. Also by Lemma 2.1, $Qf=Q((\varphi\circ q).\;\frac{g}{Qg\circ q})=\varphi . Q(\frac{g}{Qg\circ q})=\varphi$.\qed

\begin{prop}
If $\varphi\in C_c(\K)$ then there is $f\in C_c(G)$ such that $Qf=\varphi$ and $q(supp(f))\su supp(\varphi)$. Moreover if $\varphi\geq 0$ we may choose $f\geq 0$.
\end{prop}
{\bf Proof} Given $\varphi\in C_c(\K)$, there is $g\in C_c(G)_{+}$ such that $Qg=1$ on $supp(\varphi)$. put $f=(\varphi\circ q).g$ and use Lemma 2.1.\qed

\begin{prop}
If $\nu$ is a measure on $G$ satisfying
$$
\int_G f(k^{-1}xh^{-1})d\nu(x)=\Delta_K(k)\Delta_H(h)\int_G f(x)d\nu(x)\quad (k\in K,h\in H,f\in C_c(G)),
$$
then there is a unique measure $\mu=\dot\nu$ on $\K$ such that
\begin{align*}
\int_G f(x)d\nu(x) &=\int_{\K} Qf(q(x))d\dot\nu(q(x))\\
&=\int_{\K}\int_K\int_H f(kxh)dkdhd\mu(\dot{x})\quad (f\in C_c(G)).
\end{align*}
\end{prop}
{\bf Proof}
Fix $f\in C_c(G)$, then for $g\in C_c(G)$ we have
\begin{align*}
\int_G f(x)Qg(\dot x)d\nu(x) &= \int_G f(x)\int_K\int_H g(kxh)dkdhd\nu(x)\\
& =\int_K\int_H\int_G f(x)g(kxh)d\nu(x)dkdh \\
&=\int_{K}\int_H \Delta_K(k^{-1})\Delta_H(h^{-1})\int_G f(k^{-1}xh^{-1})g(x)d\nu(x)dkdh\\
&= \int_{G} g(x)\int_K \int_H f(kxh)dhdkd\nu(x)\\
&=\int_G g(x)Qf(\dot x)d\nu(x).
\end{align*}
Now using Lemma 2.3, choose $g$ such that $Qg=1$ on $q(supp(f))$, then if $Qf=0$ we have
$$
\int_G f(x)d\nu(x)=\int_G f(x)Qg(\dot x)d\nu(x)=\int_G g(x)Qf(\dot x)d\nu(x)=0,
$$
so the map $\dot\nu(Qf)=\int_G fd\nu$ is well defined on $C_c(\K)$. Given compact subset $E$ of $\K$, by Lemma 2.2 choose a compact subset $F$ of $G$ such that
$q(F)=E$ and then, by Proposition 2.1, choose $h\in C_c(G)$ supported in $F$ such that $Qh=1$ on $E$. Put
$f^{'}=(Qf\circ q)h$, then $f^{'}$ is supported in $F$, $Qf^{'}=Qf$, and
$$
|\dot\nu(Qf)|=|\dot\nu(Qf^{'})|=|\int_G f^{'} d\nu|=|\int_F f^{'} d\nu|\leq \nu(F)\|f^{'}\|_\infty \leq \nu(F)\|h\|_\infty\|Qf\|_\infty ,
$$
so $\dot\nu$ is a bounded linear functional on $C_c(\K)$ and the corresponding measure $\mu$ has the desired property. The uniqueness follows from
Proposition 2.1.\qed

\vspace{.3 cm}
The above condition leads us to the concept of {\it generalized Bruhat function} $\rho$ as described in the next proposition. We would deal with the problem of
existance of such functions later.

\begin{prop}
If $\rho:G\rightarrow\mathbb C$ is a continuous \fn , then $d\nu(x)=\rho(x)dx$ satisfies the condition of above proposition if and only if
$$
\rho(kxh)=\rho(x)\frac{\Delta_K(k)\Delta_H(h)}{\Delta_G(h)}\quad (k\in K,h\in H,x\in G,f\in C_c(G)).
$$
\end{prop}

{\bf Proof}
The condition of the above proposition for above $\nu$ could be writen as
$$
\int_G f(k^{-1}xh^{-1})\rho(x)dx=\Delta_K(k)\Delta_H(h)\int_G f(x)\rho(x)dx,
$$
so the result follows from the following calculation
$$
\int_G f(k^{-1}xh^{-1})\rho(x)dx=\Delta_G(h)\int_G f(k^{-1}x)\rho(xh)dx=\Delta_G(h)\int_G f(x)\rho(kxh)dx.\qed
$$

\vspace{.3 cm}
Now we are ready to prove a version of Weil's formula for double coset spaces.
\begin{theo}
Assume that $K$ is unimodular. There is a $G$-invariant Radon measure $\mu$ on $\K$ if and only if $\Delta_G\upharpoonleft_H=\Delta_H$. 
In this case $\mu$ is
unique up to constant factors and if suitably chosen then
\begin{align*}\int_G f(x)dx &=\int_{K\backslash G/H} Qf(q(x))d\mu(q(x))\\&=\int_{K\backslash G/H}\int_K\int_H f(kxh)dkdhd\mu(\dot{x})\quad (f\in C_c(G)).\end{align*}
\end{theo}
{\bf Proof} If $\mu$ exists then the uniqueness and above relation follow from the fact that the (left) Haar measure on $G$ (up to constant factors) and the
fact that $f\mapsto \int_{\K} Qf d\mu$ is clearly a left invariant positive linear functional on $C_c(G)$. Now again if $\mu$ exists, then for each
$k_0\in K, h_0\in H$ and $f\in C_c(G)$ we have
\begin{align*}
\Delta_{G}(h_0)\int_G f(x)dx &=\int_G f(k_0 ^{-1}xh_0 ^{-1})dx\\
&=\int_{\K}\int_K \int_H f(k_0 ^{-1}kxhh_0 ^{-1})dhdkd\mu(\dot{x})\\
&=\Delta_{H}(h_0)\int_{\K}\int_K \int_H f(kxh)dhdkd\mu(\dot{x})\\
&=\Delta_{H}(h_0)\int_G f(x)dx.
\end{align*}
Therefore, choosing $f$ appropriately, we get $\Delta_{G}(h_0) =\Delta_{H}(h_0)$.
Conversely suppose that $\Delta_G\upharpoonleft_H=\Delta_{H}$, then the result follows from Propositions 2.2 (with $\rho =1$) and 2.1.\qed

\begin{cor}
If $K$ and $H$ are compact then there is a $G$-invariant Radon measure on
$\K$.\qed
\end{cor}

\vspace{.3 cm}
When $K=\{1\}$, the above theorem gives the well-known Wiel's formula [Fl, 2.49]. When $K=H$ we get the following result.

\begin{cor}
If $H$ be a compact subgroup of $G$, then there is a $G$-invariant Radon measure $\mu$ on $G//H$ . In this case $\mu$ is unique up to constant factors and if suitably chosen then
\begin{align*}
\int_G f(x)dx=\int_{G//H}\int_H\int_H f(kxh)dkdhd\mu(\dot{x})\quad (f\in C_c(G)).\qed
\end{align*}
\end{cor}

\section{mackey machine}

There is a well developed theory of induced \repns from a closed
subgroup $H$ to $G$ using the "intermediate space" $G/H$, usually
referred to as the {\it Mackey Machine}(see for instance [Fl]).
Here we add one more ingrediant to this construction, namely a
compact subgroup $K$. Then what we get is induced \repns from $H$
to $N_G(K)$, the normalizer of $K$ in $G$, with the intermediate
space $\K$. When $K=1$, this reduces to the usual Mackey machine,
so it is reasonable to call it a generalized Mackey machine. To
get the induced representation we follow the same construction as
in the classical case [Fl].

The induction process related to $\K$ could be studied from
another point of view. This is when we consider $\K$ as a semi
convo and try to induce a \repn on $\K$ from a sub semi convo.
When $K=H$, this has been studied in a more general framework
[He], [HKK]. We briefly mention their results at the end of this
section, as we need some of them in the next section.

Let $K$ and $H$ be compact and closed subgroups of $G$,
respectively. Let $q:G\rightarrow \K$ be the quotient map.
Consider a unitary \repn $\{\sigma, \mathfrak H_\sigma\}$ of $H$.
Let
$\mathfrak F_0$ be the set of all elements $f\in C(G,\mathfrak H_\sigma)$ such that $q(supp(f))$ is compact in $\K$ and
$$
f(kxh)=\sigma(h^{-1})(f(x))\quad (k\in K,h\in H, x\in G).
$$

\begin{prop} If $\alpha\in C_c(G,\mathfrak H_\sigma)$ then
$$f_\alpha(x)=\int_K\int_H \sigma(h) \alpha(kxh) dkdh\quad (x\in G)$$
is uniformly continuous on $G$ and $\alpha\mapsto f_\alpha$ maps
$C_c(G,\mathfrak H_\sigma)$ onto $\mathfrak F_0$.
\end{prop}
{\bf Proof} Clearly $q(supp(f_\alpha))$ is included in
$q(supp(\alpha))$, and so it is compact. Also for each $h_0\in H$
we have
\begin{align*}
f_\alpha(k_0xh_0)&=\int_K\int_H \sigma(h)\alpha(kk_0xh_0h)dkdh\\
&=\int_K\int_H \sigma(h_0^{-1}h)\alpha(kxh)dkdh\\
&=\sigma(h_0^{-1})\int_K\int_H \sigma(h)\alpha(kxh)dkdh\\
&=\sigma(h_0^{-1})f_\alpha(x).
\end{align*}
Also since $K$ is compact, the function
$$\alpha^{'}(x)=\int_K\alpha(kx)dk \quad (x\in G)$$
is continuous of compact support, and so by [Fl, 6.1],
$$f_\alpha(x)=\int_H \sigma(h)\alpha^{'}(xh)dh$$
is uniformly continuous on $G$. Finally if $f\in\mathfrak F_0$,
then by Lemma 2.1 there is $\psi\in C_c(G)$ with $Q(\psi)=1$ on
$supp(f)$, and if we put $\alpha=\psi f$, then
\begin{align*}
f_\alpha(x)&=\int_K\int_H \psi(kxh)\sigma(h)f(kxh)dkdh\\
&=\int_K\int_H \psi(kxh)\sigma(h)\sigma(h^{-1})f(x)dkdh\\
&=f(x)Q(\psi)(x)=f(x).\qed
\end{align*}

\vspace{.3 cm}
Now let $N=N_G(K)=\{x\in G: x^{-1}Kx=K\}$ be the normalizer of $K$
in $G$, then $N$ acts on $\mathfrak F_0$ by left translations.
\begin{lemm} For each $x\in N$ and $f\in \mathfrak F_0$, we have
$L_x(f)\in\mathfrak F_0$ and $f\mapsto L_x(f)$ is a bijective
linear map.
\end{lemm}
{\bf Proof} Given $f\in \mathfrak F_0 $, the set
$q(supp(L_x(f)))=L_{q(x)}(q(supp(f)))$ is compact, since
$L_{q(x)}$ is continuous. Also for each $k\in K, h\in H, y\in
N,x\in G$ we have
$$
L_y(f)(kxh)=f(y^{-1}kxh)=f((y^{-1}ky)y^{-1}xh)=\sigma(h^{-1})f(y^{-1}x)=\sigma(h^{-1})L_y(f)(x),
$$
so $\mathfrak F_0 $ is stable under left translation by $G$. The
other assertions are trivial.\qed

\vspace{.3 cm}
Now let us assume that $\Delta_G \upharpoonleft_H=\Delta_H$, so that there is a $G$-invariant measure $\mu$ on $\K$.
Let $f,g\in\mathfrak F_0 $, consider
$h(x)=<f(x),g(x)>_\sigma \quad (x\in G)$, then for each $k\in K, h\in H, x\in G$
$$h(kxh)=<f(kxh),g(kxh)>_\sigma =<\sigma(h^{-1})f(x),\sigma(h^{-1})g(x)>_\sigma
=<f(x),g(x)>_\sigma,$$
so we may regard $h$ as a continuous \fn of compact support on
$\K$, and thereby we may define
$$
<f,g>=\int_{\K} <f(x),g(x)>_\sigma d\mu(\dot x),
$$
which is an inner product on $\mathfrak F_0 $. Also we have $<L_x(f),L_x(g)>=<f,g>$, for each $x\in N$.
Let $\mathfrak F$ be the Hilbert space completion of $\mathfrak
F_0$. Then for each $x\in N$, the translation operator $L_x$
extends to a unitary operator on $\mathfrak F$, and $x\mapsto
L_x(f)$ is a continuous map from $N$ to $\mathfrak F$, for each
$f\in \mathfrak F_0 $. Now as operators $L_x$ are uniformly
bounded, they are strongly continuous on all of $\mathfrak F$, and
so $x\mapsto L_x(f)$ is a continuous unitary \repn of $N$ in
$\mathfrak F$, which is denoted by $ind_H^{N_G(K)}\sigma$ and is
called the induced \repn of $\sigma$ from $H$ to $N_G(K)$.

\begin{rem}
There is a more general framework introduced by L. Pavel, which
could be used to get almost the same construction [Pl]. To show
that $\K$ fits into Pavel's framework one may use Theorem ? of
[Li] (one should put $H, G$, and $X$ of the definition in page
$434$ of [Pa] equal $H, N_G(K)$, and $G$ in our notation).
\end{rem}

In the second part of this section we briefly study those representations of $G//H$, as a hypergroup, which are canonically induced by representations 
of its subhypergroups. We quote some of the results in [HKK],[He]. As before $H$
is a closed subgroup of $G$ with the left Haar measure $\sigma$ (normalized in the compact case).

\begin{prop}

Let $\{\pi,\mathfrak H_\pi\}$ be a \repn of $G$, then

$(i)$ If $\pi(\sigma)\neq 0$, then $\pi$ gives rise to a representation $\dot\pi$ of $G//H$ given by $\dot\pi=\pi(\sigma)\pi(.)\pi(\sigma)$ on
$\pi(\sigma)(\mathfrak H_\pi)$. If $\pi$ is irreducible, then so is $\dot\pi$.

$(ii)$ $\pi(\sigma)\neq 0$ if and only if the trivial \repn $1_H$ of $H$ is a subrepresentation of $\pi\upharpoonleft_H$.\qed
\end{prop}

\vspace{.3 cm} If $(i)$ holds we say that $\pi$ is an {\it extension} of $\dot\pi$ to $G$.

\begin{prop}
Let $\{\tilde\pi, \zeta, \mathfrak H\}$ be a cyclic \repn of $G//H$ and 
$$\psi(x)=<\tilde\pi(\dot x)\zeta,\zeta> \quad (x\in G),$$ 
then the following are equivalent:

$(i)$ There is a \repn $\pi$ of $G$ such that $\tilde\pi=\dot\pi$.

$(ii)$ $\psi\in P(G)$.\qed
\end{prop}

\vspace{.3 cm}

Now let $G_1$ be a closed subgroup of $G$ containing $H$. We want to consider those \repns $\theta$ of $G_1//H$ which induce a \repn
$\theta^{'}=ind_{G_1//H}^{G//H} \theta$ of $G//H$. These are so called {\it inducible} \repns . 
Let $P:C_c(G//H)\rightarrow C_c(G_1//H)$ be defined by
$$
P(f)=(\frac{\Delta_G}{\Delta_{G_1}})^{1/2} f\upharpoonleft_{G_1//H}\quad\quad (f\in C_c(G//H)).
$$

\begin{defi}
With the above notations, $\theta$ is called inducible if for each $\zeta\in \mathfrak H_\theta, f,g\in C_c(G//H)$, we have

$(i)$ $<\theta(P(f^* *f))\zeta,\zeta>\geq 0$,

$(ii)$ $<\theta(P(g^* *f^* *f*g))\zeta,\zeta>\leq \|f\|_1^2 <\theta(P(g^* *g))\zeta,\zeta>$.
\end{defi}

\vspace{.3 cm}

In this case the induced \repn $\theta^{'}=ind_{G_1//H}^{G//H} \theta$ of $G//H$ is defined as follows: Take $V=C_c(G//H)\otimes\mathfrak H_\theta$ and define
$$
<f\otimes\zeta,g\otimes\eta>:=<\theta(P(g^* *f)\zeta,\eta>,
$$
where the right hand side is the inner product in $\mathfrak H_\theta$. Let $N=\{v\in V: <v,v>=0\}$ and identify elements $u,v\in V$ if $u-v\in N$
(take the quotient), then take the completion of this quotient with respect to the norm defined by the above inner product to get the Hilbert space
$\mathfrak H_{\theta^{'}}$ and define
$$
\theta^{'}(f)(g\otimes\zeta)=(f*g)\otimes\zeta\quad\quad (f,g\in C_c(G//H, \zeta\in \mathfrak H_\theta).
$$
This uniquely extends to a bounded linear operator on  $\mathfrak H_{\theta^{'}}$, giving a \repn of $C_c(G//H)$, which by $(ii)$ extends to a
nondegenerate \repn of $L^1(G//H)$, denoted by $\theta^{'}=ind_{G_1//H}^{G//H} \theta$.

\begin{prop} If $H$ is compact and $\dot\pi$ is a \repn of $G_1//H$ admiting an extension $\pi$ to $G_1$ which is inducible to a \repn $\pi^{'}$ of
$G$, then $\dot\pi$ is inducible to a \repn $\dot\pi^{'}$ of $G//H$ and there is an isometric embedding of $\mathfrak H_{\dot\pi^{'}}$ into
$\mathfrak H_{\pi^{'}}$. If moreover $\pi(G_1)(\pi(\sigma)(\mathfrak H_\pi))$ is total in $\mathfrak H_\pi$, then $\pi^{'}$ extends $\dot\pi^{'}$
to $G$. \qed
\end{prop}

\begin{prop} If $\lambda_{G//H}$ is the left regular \repn of $G//H$ then $\lambda_{G//H}=ind_1^{G//H} 1$ and $ind_H^G 1_H$ is the extension of
$\lambda_{G//H}$ to a \repn of $G$.\qed
\end{prop}

\section{Algebras of functions on $\K$}

In this section we study the algebras of functions on the double coset space $\K$, where $H$ and $K$ are closed subgroups of $G$. In particular, following [Fr] we
study the Fourier and Fourier-Stieltjes algebras of $\K$. The case $K=H$ is of particular interest. Recall that $G//H$ is a hypergroup with the left Haar measure
$m=\int_G \delta_{HxH} dx$ and convolution
$$
\delta_{q(x)}*\delta_{q(y)}=\int_H \delta_{q(xhy)} dh \quad (x,y\in G),
$$
and identity $q(e)$ and involution $q(x)^{-}=q(x^{-1})$, where $q:G\rightarrow G//H$ is the quotient map [Je, 8.3]. Moreover the hypergroup $G^{H}$ of orbits of the any
action of $H$ on $G$ is canonically isomorphic to $(G\times H)//(1\times H)$ [Je, 8.3]. Therefore $M(G//H)$ is a Banach $*$-algebra. In general $\K$ is only a
semi-hypergroup (semi-convo in [Je] terminology) and $M(\K)$ is only a Banach algebra.

Let $\sigma=dh$ and $\tau=dk$ be the left Haar measures on $H$ and $K$, respectively (normalised if $H$ or $K$ are compact), then put
$$
M(K:G:H)=\{\mu\in M(G): \tau *\mu *\sigma=\mu\}.
$$
The following result follows exactly like [Je, 14.2].

\begin{prop}
The Banach algebras $M(\K)$ and $M(K:G:H)$ are isometrically isomorphic. When $K=H$ theses are also isomorphic as Banach $*$-algebras.\qed
\end{prop}

Also let
$$
L^1(K:G:H)=\{f\in L^1(G): f(kxh)=f(x), \quad (a.a. x\in G, k\in K, h\in H)\},
$$
then we have
\begin{lemm}
If $K$ and $H$ are compact, then for each $f\in L^1(G)$, $fdx\in M(K:G:H)$ if and only if $f\in L^1(K:G:H)$.
\end{lemm}
{\bf Proof} Given $f\in L^1(K:G:H)$, put $\mu_f=fdx$, then for each $h\in H, k\in K$ and each Borel subset $A$ of $G$ we have
$$
\mu_f(k^{-1}Ah^{-1})=\int_{k^{-1}Ah^{-1}}f(x)dx=\int_A f(kxh)dx=\int_A f(x)dx=\mu_f(A),
$$
hence by Fubini's theorem
\begin{align*}
\tau *\mu_f *\sigma(A) &=\int_K\int_G\int_H \chi_A(kxh)dkd\mu_f(x)dh=\int_K\int_H\int_G \chi_A(kxh)dkd\mu_f(x)dh\\
& =\int_K\int_H\mu_f(k^{-1}Ah^{-1})dkdh=\mu_f(A)\int_K\int_H dkdh=\mu_f(A),
\end{align*}
that is $\mu_f\in M(K:G:H)$. Conversely if $\mu_f$ satisfies $\tau *\mu_f*\sigma=\mu_f$, then for each $g\in C_c(G)$ we have
$$
\int_K\int_G\int_H g(kxh)dkd\mu_f(x)dh=\int_G g(x)d\mu_f(x),
$$
that is
$$
\int_K\int_H\int_G [f(k^{-1}xh^{-1})-f(x)]g(x)dxdkdh=0.
$$
This implies that $f(k^{-1}xh^{-1})-f(x)=0$ for a.a $x,k,h$, that is $f\in L^1(K:G:H)$.\qed

\begin{cor}
If $K$ and $H$ are compact, then the Banach algebras $L^1(\K)$ and $L^1(K:G:H)$ are isometrically isomorphic.
When $K=H$ theses are also isomorphic as Banach $*$-algebras.\qed
\end{cor}

It is easy to see that the Banach $*$-algebra $C_b(\K)$ of bounded continuous \fns on $\K$ is isometrically isomorphic to
$$
C_b(K:G:H)=\{u\in C_b(G): u(kxh)=u(x) \quad (x\in G, k\in K, h\in H)\}.
$$
Now let $K$ be unimodular and assume that a generalized Bruhat \fn
$\rho$ exists on $G$. Let $d\nu=\rho dx$ and let $\mu=\dot\nu$ be
the corresponding quasi-invariant measure on $\K$. Define
$$
Q_\rho(f)(\dot x)=\int_K\int_H \frac{f(kxh)}{\rho(kxh)} dkdh \quad
(x\in G, f\in C_c(G).
$$
Then we have the generalized {\it Mackey-Bruhat} formula
$$
\int_{\K} Q_\rho(f)(\dot x) d\mu(\dot x)=\int_G f(x) dx.
$$
Let $\mathfrak J_\rho =\{f\in C_c(G): Q_\rho (f)=0\}$ , then

\begin{lemm}
If $K$ is unimodular, there is an isometric linear isomorphism
of Banach spaces $C_c(\K)\cong C_c(G)/\mathfrak J_\rho$ with respect to the
$\|.\|_1$-norms defined by $\mu$ and Haar measure, respectively .
\end{lemm}
{\bf Proof} It could be proved just as in the Proposition 2.1 that
the map $Q_\rho :C_c(\K)\rightarrow C_c(G)$ is surjective. Also
for each $f\in C_c(G)$ we have
\begin{align*}
\|Q_\rho(f)\|_1 & =\int_{\K} |\int_K\int_H \frac{f(kxh)}{\rho(kxh)}
dkdh |d\mu(\dot x)\\
&\leq\int_{\K} |\int_K\int_H \frac{f(kxh)}{\rho(kxh)} dkdh|d\mu(\dot x)=\int_G |f|(x)dx=\|f\|_1 .
\end{align*}
On the other hand, given $g\in C_c(\K)$, using the argument of
Proposition 2.1, we can choose $f\in C_c(G)$ such that
$Q_\rho(f)=g$ and $Q_\rho(|f|)=|g|$. Then
$$\|f\|_1=\int_G |f|dx=\int_{\K} Q_\rho(|f|)d\mu=\int_{\K}
|g|d\mu=\|g\|_1 .$$
Therefore
$$\|g\|_1=inf\{\|f\|_1 :f\in C_c(G), Q_\rho (f)=g\}=\|f+\mathfrak
J_\rho\|_1.\qed$$

\vspace{.3 cm}
Now the following result, which follows from the above lemma and [RS, Lemma 3.4.4], shows that the generalized
Mackey-Bruhat formula is valid in $L^1(G)$ also.

\begin{prop} 
If $K$ is unimodular then there is an isometric
linear isomorphism $L^1(\K, \mu)\cong L^1(G)/J_\rho$, where
$J-\rho$ is the closure of $\mathfrak J_\rho$ in $L^1(G)$.\qed
\end{prop}

Coming back to the case, $K=H$, following [Je], we denote the {\it positive definite} elements of $C_b(G//H)$ by $P(G//H)$
and its linear span by $B(G//H)$. Also, following [Fo], we consider
$$
B(G::H)=\{u\in B(G): u(kxh)=u(x)\quad (x\in G, k\in K, h\in H)\},
$$
where $B(G)$ is the Fourier-Stieltjes algebra of $G$, and denote the closure of the set of those $u\in B(G::H)$ for which $q(supp(u))\subseteq G//H$ is
compact, by $A(G::H)$. The following result is proved in [Je, 14.4].

\begin{lemm}
For each $f\in P(G)\cap C_b(G::H)$, there is $g\in P(G//H)$ such that $f=g\circ q$. \qed
\end{lemm}

\begin{prop}
The vector spaces $B(G//H)$ and $B(G::H)$ are linearly isomorphic.
\end{prop}
{\bf Proof} By above lemma, The linear map $g\mapsto g\circ q$ from $B(G//H)$ to $B(G::H)$ is surjective. It is also injective, since $q$ is surjective.\qed

\vspace{.3 cm}
The following is proved in [AM] in a more general context.

\begin{prop}
If $H$ is compact then the above linear isomorphism is norm
increasing. In particular it is an open map.\qed
\end{prop}

\begin{prop}
If $H$ is compact, then there is a surjective projection
$Q:B(G)\rightarrow B(G::H)$ with $\|Q\|\leq 1$. Also the
restriction of $Q$ on $A(G)$ is a projection onto $A(G::H)$.
\end{prop}

{\bf Proof} Put
$$
Q_1(u)(x)=\int_H u(kx)dk,\quad Q_2(u)(x)=\int_H u(xh)dh\quad (u\in
B(G)),
$$
then for $i=1,2$, $Q_i:B(G)\rightarrow B(G)$ is a continuous
projection with $\|Q_i\|\leq 1$ [Fo]. Also clearly
$Q_1Q_2=Q_2Q_1$. Hence $Q=Q_1Q_2$ is a continuous projection on
$B(G)$ with $\|Q\|\leq 1$. Now $Im(Q_2)=B(G:H)$ [Fo, 3.4], and a
similar left version holds for $Q_1$, hence $Im(Q)=B(G::H)$. The
assertion about $A(G)$ follows similarly from [Fo, 3.4].\qed

\begin{cor}
If $H$ is compact, $A(G::H)$ is a complemented subspace of both
$A(G)$ and $A(G:H)$.\qed
\end{cor}

\vspace{.3 cm} As another corollary we can prove the following
extension property. But first we need the following analogue of a
result of Herz [Hr].

\begin{lemm}
If $H$, $K$ are compact and closed subgroups of $G$ respectively,
and $H\subseteq K$. Then each $u\in A(K::H)$ extends to some $v\in
A(G)$ with the same norm.
\end{lemm}

\begin{cor}
If $H$, $K$ are compact and closed subgroups of $G$ respectively,
and $H\subseteq K$. Then each $u\in A(K::H)$ extends to some $v\in
A(G)$ with the same norm.
\end{cor}
{\bf Proof} Extend $u$ to $v_1\in A(G)$ with the same norm and put
$v=Q(v_1)$. Then $v$ extends $u$ so $\|v\|\geq\|u\|$. Also
$\|v\|\leq\|v_1\|=\|u\|$, so the result.\qed

\begin{rem}
One can show that there is a continuous surjection $Q$ as above,
even when $H$ is closed but we don't know if this is a projection.
Indeed we can argue as follows:
Recall that $f\in C_b(G)$ is {\it weakly almost periodic} if the
set $O_L(f)$ of left translates of $f$ is relatively weakly
compact. This is equivalent to its right version. We denote the
set of all such functions with $WAP(G)$. This is a
$C^*$-subalgebra of $C_b(G)$. Now let $\Psi$ be the unique
invariant mean on $WAP(G)$ and $u\in B(G)$. Let $\{x_\beta\}$ and
$\{y_\gamma\}$ be complete sets of left and right coset
representatives of $H$ in $G$. Put
$$
{}_\beta u(h)=u(x_\beta h),\quad u_\gamma(h)=u(hy_\gamma)\quad (h\in
H),
$$
and set 
$$
Q_1(u)(x)=\Psi({}_\beta u),\quad
Q_2(u)(x)=\Psi(u_\gamma)\quad (x\in x_\beta H)\quad (x\in
Hy_\gamma),
$$
then for $i=1,2$, $Q_i:B(G)\rightarrow B(G)$ is a continuous
projection with $\|Q_i\|\leq 1$. Then $Q=Q_1Q_2$ is the required
map. We are not able to show that $Q_1Q_2=Q_2Q_1$, ensuring that
$Q$ is indeed a projection.
\end{rem}

\vspace{.3 cm}

The following is proved similar to [Fo, 3.5,3.6]
\begin{prop}
If $H_1$ and $H_2$ are compact subgroups of $G$, then
$B(G::H_1)=B(G::H_2)$ if and only if $A(G::H_1)=A(G::H_2)$ if and
only if $H_1=H_2$.\qed
\end{prop}

Let $H$ be a compact subgroup of $G$. Consider
$Q:L^1(G)\rightarrow L^1(G)$ defined by
$$
Q(f)(x)=\int_H\int_H f(kxh)dkdh \quad (x\in G, f\in L^1(G)).
$$

\begin{lemm}
$Im(Q)=L^1(G::H)$ is a self-adjoint subalgebra of $L^1(G)$.
\end{lemm}

Let $C^*(G::H)$ and $VN(G::H)$ denote the norm and weak${}^*$
closure of $L^1(G::H)$ in $C^*(G)$ and $VN(G)$, respectively. The
following result is interesting as its analogue for one sided
cosets fails (c.f. remarks before Theorem 4.1 in [Fo].

\begin{prop}
If $H$ is compact,then

$(i)$ $C^*(G::H)$ is a $C^*$-subalgebra of $C^*(G)$ and $Q$
extends to a continuous projection on $C^*(G)$ with range
$C^*(G::H)$.

$(ii)$ $VN(G::H)$ is a vonNeumann subalgebra of $VN(G)$ and $Q$
extends to a continuous projection on $VN(G)$ with range
$VN(G::H)$.
\end{prop}
{\bf Proof} The first part of both assertions follow from above
lemma. The second part could be proved from the analogous result
in [Fo] and its left version using the method of Proposition
4.5.\qed

\begin{theo}
If $H$ is compact, the Banach spaces $C^*(G::H)^*$ and $B(G::H)$
are isometrically isomorphic.
\end{theo}
{\bf Proof} We have
$$
C^*(G::H)^*=Im(Q)^*\cong (C^*(G)/ker(Q))^*\cong (ker(Q))^\perp ,
$$
so we only need to check that $(ker(Q))^\perp =B(G::H)$. Let $u\in B(G::H)$ and $f\in L^1(G)$ with $Q(f)=0$. Then
\begin{align*}
<u,f> &=\int_G u(x)f(x)dx \\
&=\int_{G//H}\int_H\int_H u(kxh)f(kxh)dkdhd\mu(\dot x)\\
&= u(x)\int_{G//H}\int_H\int_H f(kxh)dkdhd\mu(\dot x)\\
&= u(x)\int_{G//H} Q(f)(\dot x)d\mu(\dot x)=0.
\end{align*}
Conversely if $u=0$ on $ker(Q)$, then for each $k_0, h_0\in H$ let
$g(x)=f(k_0^{-1}xh_0^{-1})-f(x)$, then clearly $g\in ker(Q)$ so we
have
$$
\int_G (u(k_0 xh_0)-u(x))f(x)dx=\int_G u(x)g(x)dx=0.
$$
This being true for each $f$, we have $u(k_0 xh_0)-u(x)=0$, for
each $x\in G$, that is $u\in B(G::H)$. \qed

\vspace{.3 cm} Similarly we have
\begin{theo}
If $H$ is compact, the Banach spaces $VN(G::H)_{*}$ and $A(G::H)$
are isometrically isomorphic.\qed
\end{theo}

\begin{prop}
Let $H$ be compact, then the following are equivalent.

$(i)$ $G$ is amenable.

$(ii)$ $A(G::H)$ has a bounded approximate identity consisting of
functions of compact support.

$(iii)$ $A(G::H)$ is weakly factorized.
\end{prop}
{\bf Proof} $(i)\Rightarrow (ii)$ Given $E\subseteq G$ compact, let
$F=HEH$, then $F\subseteq G$ is compact and so by Reiter's
condition, for each $\epsilon>0$, there is an element
$u=u_{E,\epsilon}\in A(G)$ with compact support such that
$u=(1+\epsilon)^{-1}$ on $F$ and $\|u\|\leq 1$. Let
$v=v_{E,\epsilon}=Q(u)\in A(G::H)$, then clearly
$v=(1+\epsilon)^{-1}$ on $E$ and $supp(v)\subseteq H(supp(u))H$,
so $v$ is of compact support. Let $w\in A(G::H)$ be any element of
compact support, and let $E=supp(w)$, then for $v=v_{E,\epsilon}$
we have $vw=(1+\epsilon)^{-1}w$, so
$\|vw-w\|=\|w\|.|(1+\epsilon)^{-1}-1|<\epsilon\|w\|$. Since
elements of $A(G::H)$ with compact support are dense in $A(G::H)$,
this proves $(ii)$.

$(ii)\Rightarrow (iii)$ Cohen's factorization theorem.

$(iii)\Rightarrow (i)$ Let $F\subseteq G$ be compact, then by weak
factorization property, there is an element $u=u_F\in A(G::H)$
with $u>1$ on $F$ and $\|u\|\leq M$, for some constant $M$
independent of $F$ [Fo]. Take any $f\in C_c(G)_{+}$ and let
$F=supp(f)$ and observe that
$$
<u_F,f>=\int_G u_F(x)f(x)dx \geq \int_G f(x)dx=\|f\|_1 .
$$
The rest of the proof is standard and could be finished as in
[Fo,4.2].\qed

\begin{theo} If $G$ is amenable and $H$ is compact, then the
multiplier algebra $M(A(G::H))$ is isometrically isomorphic to
$B(G::H)$.
\end{theo}
{\bf Proof} By above theorem, $A(G::H)$ has a \bai $\{u_\alpha\}$
(bounded by $1$). Let $M=M(A(G::H))$ and $\|.\|_M$ denote the
multiplier norm. It is clear that $B(G::H))\subseteq M\subseteq
C_b(G)$. Take $u\in M$, then $uu_\alpha\in A(G::H)$ and
$\|uu_\alpha\|\leq \|u\|_M \|u_\alpha\|\leq \|u\|_M$, for each
$\alpha$. Now $\{uu_\alpha\}$ converges to $u$ in the norm of
$A(G)$ and so pointwise. Hence $u\in B(G::H)$ and
$\|u\|\leq\|u\|_M$. The inequality in other direction now follows
from $B(G::H)\subseteq M$. \qed

\end{document}